\documentclass[twoside]{article}
\usepackage[a4paper]{geometry}
\usepackage[latin1]{inputenc} 
\usepackage[T1]{fontenc} 
\usepackage{RR}
\usepackage{hyperref}
\usepackage{algorithmic,algorithm}
\usepackage{amsmath,amssymb} 
\usepackage{graphicx}
\usepackage{booktabs}
\newcommand{\SeaSpace}{\Omega}
\newcommand{\Nbh}{\mathcal{N}}

\newcommand{\Real}{\mathop{\rm I\kern-.2emR}}
\newcommand{\V}{{\cal N}}
\RRNo{8215}
\RRdate{Januray 2013}
\RRauthor{
Marie-El\'eonore Marmion \thanks{Universit\'e Libre de Bruxelles, Belgium \url{(marie-eleonore.marmion@ulb.ac.be)}}
\and Aymeric Blot \thanks{ENS Cachan/Bretagne, Universit\'e Rennes 1, France \url{(blot.aymeric@gmail.com)}}
  \and
Laetitia Jourdan \thanks{Universit\'e Lille 1 - INRIA Lille Nord-Europe, France \url{(laetitia.jourdan@lifl.fr)}} 
\and
Clarisse Dhaenens \thanks{Universit\'e Lille 1 - INRIA Lille Nord-Europe, France \url{(clarisse.dhaenens@lifl.fr)}} 
}
\authorhead{Marmion \& others}
\RRtitle{Etude de la Neutralit\'e dans le probl\`eme de Coloration de graphe}
\RRetitle{Neutrality in the Graph Coloring Problem}
\titlehead{Neutrality in the Graph Coloring Problem}

\RRresume{
Le probl\`eme de coloration de graphes est souvent \'etudi\'e dans la litt\'erature et est r\'eguli\`erement not\'e que beaucoup de solutions voisines ont la m\^eme qualit\'e.
Cependant, \'a notre connaissance, aucune analyse approfondie de cette neutralit\'e n'a jamais \'et\'e men\'ee dans la litt\'erature.
Dans cet article, nous quantifions la neutralit\'e de certaines instances difficiles du probl\`eme de coloration de graphe.
Il est int\'eressant de d\'etecter cette propri\'et\'e de neutralit\'e car elle influe le processus de recherche.
En effet, les optima locaux peuvent appartenir \'a des plateaux qui repr\'esentent un obstacle pour les m\'ethodes de recherche locale.
Dans ce travail, nous visons \'egalement \'a montrer l'int\'er\^et d'exploiter la neutralit\'e au cours de la recherche. Par cons\'equent, les performances de NILS, une recherche locale g\'en\'erique et d\'edi\'ee aux probl\`emes neutres, sont \'etudi\'es sur plusieurs instances difficiles.
}
\RRabstract{
The graph coloring problem is often investigated in the literature.
Many insights about many neighboring solutions with the same fitness value are raised but as far as we know, no deep analysis of this neutrality has ever been conducted in the literature.
In this paper, we quantify the neutrality of some hard instances of the graph coloring problem.
This neutrality property has to be detected as it impacts the search process.
Indeed, local optima may belong to plateaus that represents a barrier for local search methods.
In this work, we also aim at pointing out the interest of exploiting neutrality during the search. Therefore, a generic local search dedicated to neutral problems, NILS, is performed on several hard instances.  }
\RRmotcle{Graph Coloring Problem, Fitness Landscape, Neutrality}
\RRkeyword{Graph Coloring Problem, Fitness Landscape, Neutrality}
\RRprojet{Dolphin}
\RCLille 

\begin{document}
\makeRR   
\tableofcontents
\section{Introduction}

The Graph Coloring Problem (GCP) is a well known combinatorial problem defined as follows:
considering a graph $G=(V,E)$ where $V$ is the set of $n$ vertices (nodes) and 
$E$ is the set of edges ($E \in V \times V$),
the goal is to find the minimum number $k$, such that, given the mapping $\Psi : V \mapsto 1,2...k$,
for all edge $(u,v) \in E$, $\Psi(u)\neq\Psi(v)$.
Such $k$ is called the chromatic number of $G$ and is denoted by $\chi$. \\

The GCP is a $\mathcal{NP}-$hard problem~\cite{garey:1990} widely studied in the literature.
Then, for larger instances, approximate algorithms are used.
For example, a lot of local search methods have been applied or proposed to solve the GCP~\cite{avanthay:2003,caramia:2006,galinier:2006,hertz:1987,porumbel:2010}.
Moreover, the GCP has also been tackled by evolutionary strategies~\cite{porumbel:2010_2}.
Indeed, the most efficient metaheuristic schemes have been adapted to the GCP, including specific encodings and mechanisms, in order to give better performance. All these adaptations require a very good knowledge of the problem and a long time of experimental analysis to tune the best parameters. \\

Another way to design efficient algorithms is to analyze the problem structure.
Therefore, let us, first, remember the basic notions to define the structure of a combinatorial optimization problem.
The \emph{search space} $\SeaSpace$ is the set of admissible solutions
and $f: \SeaSpace \longrightarrow \Real$ is a \textit{fitness function} that assigns a quality to each solution s $\in \SeaSpace$.
A \emph{neighborhood structure} is a mapping function $\Nbh: \SeaSpace \rightarrow 2^\SeaSpace$ 
that assigns a set of solutions $\Nbh (s) \subset \SeaSpace$ to any feasible solution $s \in \SeaSpace$.
$\Nbh(s)$ is called the \emph{neighborhood} of $s$, 
and a solution $s' \in \Nbh(s)$ is called a \emph{neighbor} of $s$.
A \emph{fitness landscape} \cite{stadler:1996,wright:1932} can be defined by a triplet $(\SeaSpace, \Nbh, f)$. Landscape may be a way to describe the problem structure.
Many authors investigated the landscape of different combinatorial optimization problems~\cite{bouziri:2011,marmion@jmma:2011,merz:1997,merz:2000_2}.
This landscape analysis aims at understanding better the characteristics of the problems
and then, designing efficient algorithms. 
For example, neutrality appears when neighboring solutions have the same fitness value.
Then, neutrality is a characteristic of the landscape~\cite{reidys:2001}
and has been analyzed on different problems such as the permutation flowshop scheduling problem with the makespan criterion~\cite{marmion@lion:2011} and the NKq-landscape problem~\cite{verel:2004}
in order to design effective local search to solve problem under neutrality. \\

The GCP is often investigated in the literature.
Many insights about the neutrality of this problem are raised
when considering the number of edges with the same color at both ends.
But, as far as we know, no deep analysis has ever been conducted in the literature.
In this paper, we aim at investigating if the GCP may be considered as a neutral problem and if the neutrality may be exploited to solve the GCP.
Therefore, in Section~2, the $k$-GCP is defined and the modeling used in this article is given.
Section~3 gives measures to analyze the neutrality of a combinatorial optimization problem and
gives the results on the neutrality study of the GCP instances.
Then, in Section~4, the benefit of exploiting the neutrality when solving the GCP is studied.
Section~5 gives the conclusions of the presented work and future research interests.

\section{The Graph Coloring Problem}

The GCP is a problem widely studied in literature.
Indeed, many real problems are modeled and solved using graph coloring. 
Among them, many applications are found such as frequency assignment to antennas~\cite{maniezzo:2000}, schedule design~\cite{zufferey:2008} or register allocation~\cite{briggs:1994}. \\

As the GCP is a NP-hard problem, several heuristics and metaheuristics have been proposed to solve the large instances. These approaches can be classified in three main solution approaches~\cite{porumbel:2010}: (i) sequential construction heuristics like Dsatur~\cite{brelaz:1979} that are fast but not really efficient, (ii) local search algorithms and (iii) evolutionary hybrid or population-based algorithms.  
For local search algorithms, several methods using problem-specific heuristics have been proposed~\cite{galinier:2006} such as Tabu Search~\cite{hertz:1987,fleurent:1996}, Simulated Annealing~\cite{chams:1987}, Iterated Local Search~\cite{paquete:2002}, VNS~\cite{tricjk:2007}, etc. In particular, the CHECKCOL algorithm~\cite{caramia:2006}
found new chromatic numbers for four different DIMACS instances. Recently, Porumbel \textit{et al.}~\cite{porumbel:2010} have also proposed two tabu search algorithms (TS-Div and TS-Int) based on the TABUCOL~\cite{hertz:1987} heuristic leading to good performance. Moreover, Avanthay \emph{et al.}~\cite{avanthay:2003} have adapted and applied a variable neighborhood search to solve the GCP efficiently.
In addition, evolution strategies (ES) are one example of population based approaches that have been proposed to solve the GCP. In particular a very competitive algorithm on DIMACS instances is an hybrid algorithm which combines a tabu search with a genetic algorithm~\cite{porumbel:2010_2}. \\

In the following, we specify exactly the problem considered and present the instances used in the experiments.
 
\subsection{Problem Definition and Representation used}
The GCP consists in finding the minimal number of colors $\chi$, called the chromatic number, 
that leads to a legal coloring of a graph. 
The $k$-GCP is a related problem, that deals with the existence of a legal coloring using $k$ colors ($k$-coloring). 
Since it is easier to find a coloring with $k > \chi$, 
the following strategy is often used to solve GCP via $k$-GCP:
\begin{itemize}
\item[(i)] Generate an initial legal $k$-coloring ($k > \chi$)
\item[(ii)] Set $k = k - 1$
\item[(iii)] Solve the $k$-GCP; if a legal $k$-coloring is found
then go step (ii), else return~$k$
\end{itemize}
In this study, we are interesting in analyzing the neutrality of the $\chi$-GCP problem, that is to say, the most difficult $k$-GCP problem since the best known number of colors is $\chi$.
Then, characterizing the structure of the $\chi$-GCP should help to find a solution 
that represents a legal coloring with $\chi$ colors.\\

In order to lead a landscape analysis of the problem we need to define the three elements  $(\SeaSpace, \Nbh, f)$, where $\SeaSpace$ is the search space and depends on the representation, $\Nbh$ is a neighborhood induced by the neighborhood operator and $f$ is the objective function. In the literature, several representations, neighborhood and objective functions are found to deal with the GCP.\\
In the case of $\chi$-GCP, we propose to adopt the following representation:
A solution is represented as a vector of colors, $s=[c(1),\ldots,c(i),\ldots,c(n)]$
where $n$ is the number of nodes and $c(i)$ is the color associated to node $i$.
To ensure that each solution has one unique representation, colors
are enumerated by their order of arrival. Therefore, the colorings
that differ only by a permutation of colors have the same representation. The search space  $\SeaSpace$ is then defined by all the possible vectors.
The \emph{$1$-move} operator defines the neighorhood relation. This operator changes the color of one node. With $n$ nodes and $k$ colors, the neighborhood size of a solution is bounded by $n \times (k-1)$.
The objective function $f$ aims at evaluating how far the solution is from a legal coloring. In this work, we choose to associate to a solution the number of conflicts. Hence, the fitness value of a solution is equal to the number of edges with the two endpoints of same color:
$$f : s \to \sum_i|E_i|$$ where
$E_i$ is the set of edges with both endpoints of color $i$ in the solution $s$.

\footnotetext[1]{\url{http://dimacs.rutgers.edu/Challenges/}}

\subsection{Benchmark Problems}

In this study we focus on literature instances known to have ``difficult upper bound'', 
that is to say that a minimal legal coloring is hard to obtain. 
Those instances are extracted from the DIMACS Computational Challenge on ``Graph Colouring and its Generalisations''\footnotemark[1].
There are four classes of instances, according to the type of generation.
\begin{itemize}
\item {dsjc$X$.$Y$} are graphs with $X$ vertices, where $Y$ is the probability that two vertices are connected by an edge.
\item {dsjr$X$.$Y$ and r$X$.$Y$} are graphs with $X$
vertices. Two vertices are connected if their distance is less than
$Y$. A suffix ``c'' denotes the complementary of the graph.
\item {flat$X$\_$Y$} are graphs with $X$ vertices, based on
an initial $Y$-classes partitioning. Finding the best legal coloring
is equivalent to restoring the initial partitioning.
 \item {le$X$\_$Y$} are graphs with $X$ vertices, based on a clique of size $Y$. 
For those graphs, $\chi = Y$.

\end{itemize}

In the next section, a study under the neutrality point of view of the structure of the
DIMACS instances for the $\chi$-GCP is leaded.

\section{Neutrality in GCP}

In this work, we are interested in the neutrality property.
First, definitions and measures are given in order to characterize the neutrality of
a combinatorial optimization problem such as the graph coloring problem.
Then, experiments are led on the GCP to analyze the neutrality of the different instances.

\subsection{Measures to characterize the neutrality}

A \emph{neutral neighbor} of a solution is a neighboring solution having the same fitness value.
The set of neutral neighbors of a solution $s \in \SeaSpace$ is then $\Nbh _n(s) = \{ s' \in \Nbh (s) ~|~ f(s^\prime) = f(s) \}$.
The \emph{neutral degree} of a given solution is the number of neutral solutions in its neighborhood.
A fitness landscape is said to be \emph{neutral} if there are ``many'' solutions with a high neutral degree. 
A \emph{neutral fitness landscape} can be pictured by a landscape with many plateaus.
The average or the distribution of neutral degrees over the landscape may be used to qualify the level of neutrality of a problem instance. 
This measure plays an important role in the dynamics of local search algorithms \cite{verel:2007,wilke:2001}.\\

\begin{figure}[ht]
\centering
\includegraphics[height=2cm]{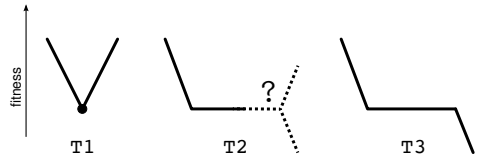}
\caption{Typologies of plateaus.}
\label{plateauGeometry}
\end{figure}

In the case a problem gets the neutrality property, Marmion \emph{et al.}~\cite{marmion@lion:2011}
suggested to characterize the plateaus found from the local optima. 
Indeed, theses plateaus might trap a local search even though all the solutions belonging to such a plateau
are not necessarily local optima.
Thus, we define by the term \textit{portal} a solution in a plateau of a local optimum, having at least one neighbor with a better fitness value.
When the fitness landscape is neutral, an important characteristic of the landscape can be described by its plateaus, that may be sampled by \textit{neutral random walks}.
A neutral random walk is a sequence of solutions where the solution $s_{i+1}$ is randomly chosen in the neutral neighborhood of the solution $s_{i}$.
The plateaus of the local optima sampled by a neutral random walk were classified in a three-class topology (see Figure~\ref{plateauGeometry}): 
(T1) the local optimum is the single solution of the plateau, {\itshape i.e.} it has no neutral neighbor; 
(T2) no neighbor with a better fitness value was met for any solutions of the plateau encountered along the neutral walk and; 
(T3) a portal has been identified on the plateau, {\itshape i.e.} at least, one solution of the plateau has an improving neighbor.\\

Let us give the measure to characterize the neutrality through 
the study of the plateaus of the local optima.
An important information is given by the \textit{autocorrelation of neutral degree} along a neutral random walk \cite{bastolla:2003} as it measures the correlation of the plateau structure. The autocorrelation function $\rho(k)$ \cite{weinberger:1990} is the correlation coefficient of the neutral degree between the solutions $s_{i}$ and $s_{i+k}$ of the neutral walk. 
If the first correlation coefficient $\rho(1)$ is close to 1, the neutral degree variation between neighbors is low, and so plateaus may be considered as structured graph.
As neutral walks are random, they may loop on a subset of solutions.
To attest that a neutral walk correctly describes a plateau, 
the number of solutions that are encountered at least twice during the sampling is computed.
Moreover, in case of a T3 plateau, the position of portals is an interesting information.
Indeed, the number of solutions visited before finding a portal during a neutral random walk is a good indicator of the probability to find an improving solution.\\

Characterizing the plateaus of the local optima aims at understanding their role 
during the search process.
Then, there exists a cost/quality trade-off between the number of solutions visited to find a portal
and the number of solutions visited to find a (new) local optimum starting with a new solution, called the step length.
This trade-off depends on the number of T3 plateaus relative to the number of T2 plateaus
sampled for a same instance.
This trade-off consists in analyzing the fitness landscape from the local search dynamics point of view. Indeed, when a local optimum is found, the measures should help to answer
if it is faster to restart from another solution in order to find a new better local optimum
or, to exploit neutrality of the landscape to move on plateaus in order to meet a portal and then accept an improving solution.

\subsection{Experimental setup}

The average neutral degree is estimated from 30 random solutions uniformly generated
and the average neutral degree of local optima is estimated from 30 local optima
found by a steepest descent starting from random solutions.
The ratio of the neutral degree is the neutral degree over the size of the neighborhood.
This measure is also computed for the random solutions and the local optima
as it makes the comparison between different instances easier.

The step lengths of the descents to find the local optima have been recorded. 
The maximal step length value, $L_{max}$, is used to define the length of the neutral walks P$_{max}$.
Thus, the lengths are comparable to analyze the trade-off between 
the number of solutions visited to find a portal
and the number of solutions visited to find a (new) local optimum.
Then, 30 neutral walks are run independently from the 30 different local optima.
The neighborhood of each solution from the neutral walk is visited entirely
in order to compute the numbers of improving, neutral and worsening neighbors.
Thus, a given instance is characterized by 30 plateaus.

\subsection{Experimental results}

\begin{table}[t!]
  \centering
  \tabcolsep=2.5pt
   \caption{Average neutral degree and the corresponding ratio
   for the random solutions and the local optima.}
  \begin{tabular}{c c ccc c cc c cc}
  \\
  \toprule
   &\multicolumn{5}{c}{}&\multicolumn{5}{c}{Neutral Degree} \\
   \cmidrule{7-11}
    & &\multicolumn{3}{c}{Data}& &\multicolumn{2}{c}{Random}& &\multicolumn{2}{c}{Local optima}\\
    \cmidrule{3-5} \cmidrule{7-8} \cmidrule{10-11}
       Instances  & $~~~$ &  $V$ &  $\chi$ &  $\mid$nbh$\mid$&  $~~~$ & nd & ratio  &   & nd & ratio\\
    \cmidrule{1-11} 
    dsjc250.5    &   & 250& 28& 6750&      &      858&12.7\%& & 83.8&1.2\%\\
    dsjc500.1    &   & 500& 12& 5500&      &      800&14.5\%& &  144&2.6\%\\
    dsjc500.5    &   & 500& 48&23500&      &     2910&12.4\%& &  176&0.7\%\\
    dsjc500.9    &   & 500&126&62500&      &     9320&14.9\%& &  384&0.6\%\\
    dsjc1000.1   &   &1000& 20&19000&      &     2440&12.8\%& &  290&1.5\%\\
    & & \\
   
    r250.5        &  & 250& 65&16000&       &    3470&24.2\%& & 1090&9.7\%\\
    dsjr500.5     &  & 500&122&60500&       &   12800&21.1\%& &  356&5.9\%\\
    dsjr500.1c    &  & 500& 84&41500&       &    4780&11.5\%& &  121&0.3\%\\
    dsjr1000.1c   &  &1000& 98&97000&       &    8600& 8.9\%& &  188&0.2\%\\
 & & \\
 
    flat300\_28\_0 & & 300& 28& 8100&        &   1010&12.5\%& & 90.5&1.1\%\\
    flat1000\_50\_0 & &1000& 50&49000&       &    4400& 9.0\%& &   146&0.3\%\\
     & & \\
     
    le450\_25c   &   & 450& 25&10800&      &     1910&17.7\%&  & 552&5.1\%\\
    le450\_25d   &   & 450& 25&10800&      &     1900&17.6\%&  & 496&4.6\% \\
   \bottomrule 
  \end{tabular}
  \label{GCP_neutralD}
\end{table}

\begin{figure}[h!]
\centering
\includegraphics[scale=0.50]{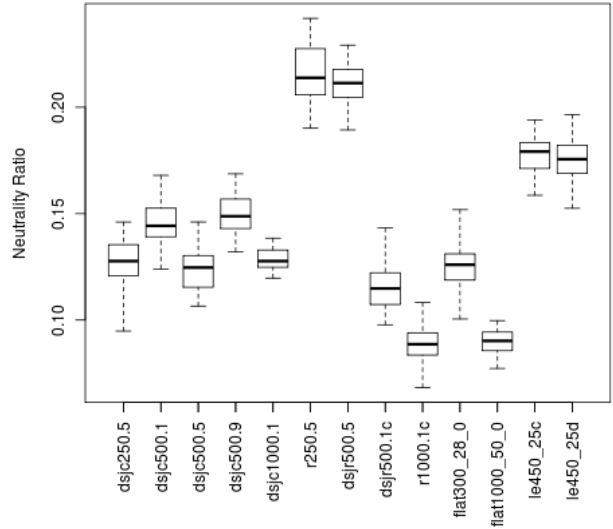}
\caption{Boxplots of the neutral degree ratio for a sampling of 30 random solutions.}
\label{GCP_neutralD_boxplot_rand}
\end{figure}

Table~\ref{GCP_neutralD} gives the average neutral degree and the corresponding ratio
for random solutions and local optima on the GCP instances.
For each instance, the number of nodes $V$, the chromatic number $\chi$ and the size of the neighborhood $\mid$nbh$\mid$ are also given. This table first shows that, the ratios for random solutions are quite high (up to 24.2\% for the instance r250.5). That indicates the neutrality property of the considered instances.
Figure~\ref{GCP_neutralD_boxplot_rand} shows the boxplots of the neutral degree ratio of the random solution.
The distribution of the ratios points out that, for an instance, the neutral degree of the solutions do not vary. 
Then, the neutrality characterises the problem in general.
The landscape may have a lot of flat parts.
The second observation on the results of Table~\ref{GCP_neutralD} is that the ratios for local optima are smaller than the ones of the random solutions. Hence, depending on the instances, the number of neutral neighbors in the neighborhood of a local optimum may be important or not.
Some instances present a high neutral degree (such as 9.7\% for the highest, or around 5\% for others) 
while some others have a much smaller neutral degree (down to 0.2\% for the instance dsjr1000.1c).
These results confirm the $a~priori$ fact that neutrality is a strong property in the graph coloring problem that may be used in local search strategies to be more efficient. \\

Therefore, the plateaus of the local optima have to be analyzed in details
when the average neutral degree of the local optima is significant. This analysis will be achieved using the \textit{autocorrelation of neutral degree} along a neutral random walk as proposed before.
In the following, only the instances where the average neutral degree of the local optima is
higher than 1\% are considered. Moreover, since the random neutral walks explore the neighborhood entirely, the computation time can be very high. Then, only the instances where the neighborhood size is lower than 16,000 solutions are analyzed. \\

\begin{table}[h!]
  \centering
    \caption{The first value of the autocorrelation of the neutral degree }
  \begin{tabular}{c c c}
  \\
  \toprule
   Instances  & $~~~$   &  $\rho(1)$ \\
   \cmidrule{1-3} \\
    dsjc250.5      &  & 0.74\\
    dsjc500.1      &  & 0.69\\
    r250.5         &  & 0.90\\
    flat300\_28\_0 &  & 0.74\\
    le\_450\_25c   &  & 0.78\\
    le\_450\_25d   &  & 0.77\\
    \bottomrule \\
  \end{tabular}
  \label{GCP_rho1}
\end{table}

First, the first value of the autocorrelation $\rho (1)$ of each neutral walk has to be checked 
in order to verify that it is reliable to characterize the structure of the plateau.
Table~\ref{GCP_rho1} gives the value of $\rho(1)$ for all considered instances.
Then, $\rho(1)$ values are higher than 0.69 
that ensures the reliability of the following results.\\

\begin{table}[h!]
  \centering
    \caption{Types of plateaus encountered by the 30 neutral walks}
  \begin{tabular}{c c ccc}
  \\
  \toprule
    &  &\multicolumn{3}{c}{Plateaus}\\
    \cmidrule{3-5}
   Instances  & $~~~$   &  T1 & T2 & T3 \\
   \cmidrule{1-5} \\
    dsjc250.5      &  & 0& 0&30\\
    dsjc500.1      &  & 0& 0&30\\
    r250.5         &  & 0& 1&29\\
    flat300\_28\_0 &  & 0& 0&30\\
    le\_450\_25c   &  & 0& 0&30\\
    le\_450\_25d   &  & 0& 0&30\\
    \bottomrule \\
  \end{tabular}
  \label{GCP_typology}
\end{table}

Table~\ref{GCP_typology} presents the number of each type of plateaus found by the 30 neutral walks.
No degenerated plateau, with a single solution, is found.
It means that the 30 local optima found have, at least, one neutral neighbor.
Only plateaus with at least a portal have been explored by the 30 neutral walks,
except for the r250.5 instance, where one neutral walk did not meet a portal.
But, it is still possible that the sampled plateau gets one. 
Moreover, no solutions is visited twice by a neutral walk.
Consequently, the neutral walks do not loop and correctly sample the plateaus. \\

These results are promising in order to exploit the neutrality in local search algorithms.
At this step, we must verify that portals may be found quickly by the neutral walks in order to consider the plateau as a new start of the search.\\

Table~\ref{GCP_portals} gives the statistics of the distance from the
initial local optimum to the closest portal found \textit{i.e.} the number of solutions visited on a plateau ($nbS$) before meeting a portal. The statistics of the step lengths ($L$) are also given to study the cost/quality trade-off.
Clearly, it is very quick to meet a portal even randomly. 
Indeed, it is necessary to visit only 1 or 2 new solutions on the plateau to find a portal to escape. 
Let us remark that reaching a portal may be quick, 
but the difficulty is to identify a solution as a portal. 
However, compared to the steps lengths, it seems to be more interesting to continue the search process by moving on a plateau to find a portal than to restart the search process from a new random solution. \\

\begin{table}[h!]
  \centering
  \caption{Cost/quality trade-off - The number of T3 plateaus is
    remembered to ensure the reliability of the statistics. nbS gives
    the number of visited solutions before finding a portal and LgM
    gives the number of visited solutions needed to find a local
    optimum starting from a random solution of the search space using a steepest descent.}
  \begin{tabular}{ccccccccccccc}
  \\
  \toprule
      & &T3& &\multicolumn{4}{c}{$nbS$}& &\multicolumn{4}{c}{$L$}\\
      \cmidrule{3-3}  \cmidrule{5-8}  \cmidrule{10-13} 
      Instances& $~~~$       &  nb& $~~~$ & Min&Med&Mean&Max& $~~~$ & Min&Med&Mean&Max\\
    \cmidrule{1-13} \\
    dsjc250.5     & & 30& & 1&1&1.7& 6&   &     266&301&301&323\\
    dsjc500.1     & & 30& & 1&2&  2& 5&   &     501&530&532&596\\
    r250.5        & & 29& & 1&2&3.3&17&   &     125&148&148&164\\
    flat300\_28\_0& & 30& & 1&1&1.7& 6&   &     344&388&385&406\\
    le\_450\_25c  & & 30& & 1&1&2.3& 9&    &    362&396&399&424\\
    le\_450\_25d  & & 30& & 1&2&2.5&11&   &     370&399&400&428\\
    \bottomrule
  \end{tabular}
  \label{GCP_portals}
\end{table}

Experiments confirm the intuition that the GCP, associated to the triplet  $(\SeaSpace, \Nbh, f)$ defined in section 2, presents the neutrality property.
Moreover, local optima have neighboring solutions with the same fitness value.
It leads to plateaus that may disturb the progress of the search process.
The analysis of the plateaus of the local optima shows that portals, solutions of the plateau with at least one improving neighbor, are quick to reach with a random neutral walk.
It assumes that exploiting neutrality in the search process may help to find better solutions. 
The following section provides insight into how to exploit neutrality to solve the GCP.

\section{Influence of Neutrality on Local Search Performance}

In the section below, the neutrality property of some hard instances of
the graph coloring problem has been highlighted.
Indeed, for these instances, the local optima have a significant number of neutral neighbors
that leads to plateaus trapping the search process. 
Here, we propose to study the interest of exploiting the neutrality by moving on the plateaus
when a local optima is found.
\emph{NILS} is an iterated local search based on neutrality~\cite{marmion@evocop:2011}, 
which has been efficiently applied to another neutral problem, 
the permutation flowshop scheduling problem with the makespan criterion.
Thus, NILS is a generic local search that benefits from the neutrality of the problem.
In this section, NILS is run on the graph coloring problem to point out
the interest of exploiting the neutrality in a local search.

\subsection{NILS algorithm}

The Neutrality-based Iterated Local Search (NILS) is an algorithm designed to exploit the plateaus
of the local optima.
This algorithm iterates a steepest descent, called \emph{First Improvement Hill Climber} (FIHC),
and  a perturbation step to escape when the search is blocked on a local optimum. 
The used FIHC replaces the current solution by the first improving neighbor found.
The neighborhood of a solution is evaluated in a random order, and each neighbor is evaluated only once. \\

There are two ways to escape from the plateau of a local optima in a neutral
fitness landscape: either performing neutral moves until finding a portal, 
or performing a kick move which is a ``large step'' move.
A neutral move should be applied when it is assumed that the \emph{exploitation} of the neutral properties helps to find a better   solution. 
On the contrary, a kick move should be applied when it is assumed that portals are hard to find and therefore \emph{exploration} of another part of the search space is more promising.
The Neutral Walk-Based Perturbation of NILS deals with \emph{exploitation} and \emph{exploration} by setting a maximum number of steps ($MNS$) allowed in a neutral walk  (see Algorithm~\ref{algoPerturb}). 
If no improving neighbor has been found until then, the solution is kicked. 
Otherwise, the neutral walk is stopped and a local search is performed from the improving neighbor.
As for FIHC, the neighborhood is evaluated in a random order, and each neighbor is evaluated only once.
In the following, NILS is performed on the graph coloring problem in order to emphasize the benefit of exploiting the neutrality of this problem.

\begin{algorithm}[t]
\caption{Neutral Walk-based Perturbation (NWP)}
\label{algoPerturb}
\begin{algorithmic}
\STATE step $\leftarrow$ 0, better $\leftarrow$ \FALSE
\WHILE{step $<MNS$ \AND \NOT better \AND $|\V_n(s)|>0$}
	\STATE choose $s^{\prime} \in \V (s)$ such that $f(s^{\prime}) \leq f(s)$
	\IF{$f(s^{\prime}) < f(s)$}
		\STATE better $\leftarrow$ \TRUE
	\ENDIF
	\STATE $s \leftarrow s^{\prime}$ 
	\STATE step $\leftarrow$ step$+1$
\ENDWHILE
\IF{\NOT better}
	\STATE $s \leftarrow$ \textbf{kick}($s$)
\ENDIF
\end{algorithmic}
\end{algorithm}

\subsection{Experimental setup}
Four instances, one of each type, have been selected to perform NILS: 
dsjc250.5, r250.5, flat\_300\_28\_0 and le\_450\_25\_c.
The landscape analysis of theses instances has shown that the plateaus of the local optima get portals that lead to improving solutions. \\

The $MNS$ value is the single parameter of NILS algorithm.
It controls the exploitation of the plateaus by allowing a maximal number of neutral moves on each plateau.
In this paper, the aim is to point out the efficiency of using neutrality to solve the graph coloring problem.
Then, several $MNS$ values were tested in order to analyze the trade-off between exploiting the plateau
and exploring an other part of the search space.
$MNS$ values were set to 1, 2 and 5 times the size of the neighborhood.
Moreover, $MNS$ values set to 0 was tested as it corresponds to a classical ILS that restarts from a new part of the search space.
This classical ILS and NILS with a positive (strict) $MNS$ value share a similar behavior
and differ only during the perturbation step when NILS continues the search by moving neutrally
when a local optimum is found.
The comparison of the performance of NILS and this classical ILS aims at
pointing out the benefit of exploiting neutrality during the search.
For each configuration of NILS and the ILS, 30 runs were performed.
The stopping criterion was set to $2 \times 10^7$ evaluations.

\subsection{Experimental results}

\begin{figure}[ht!]
\centering
  \begin{tabular}{cc}
  \includegraphics[scale=0.35]{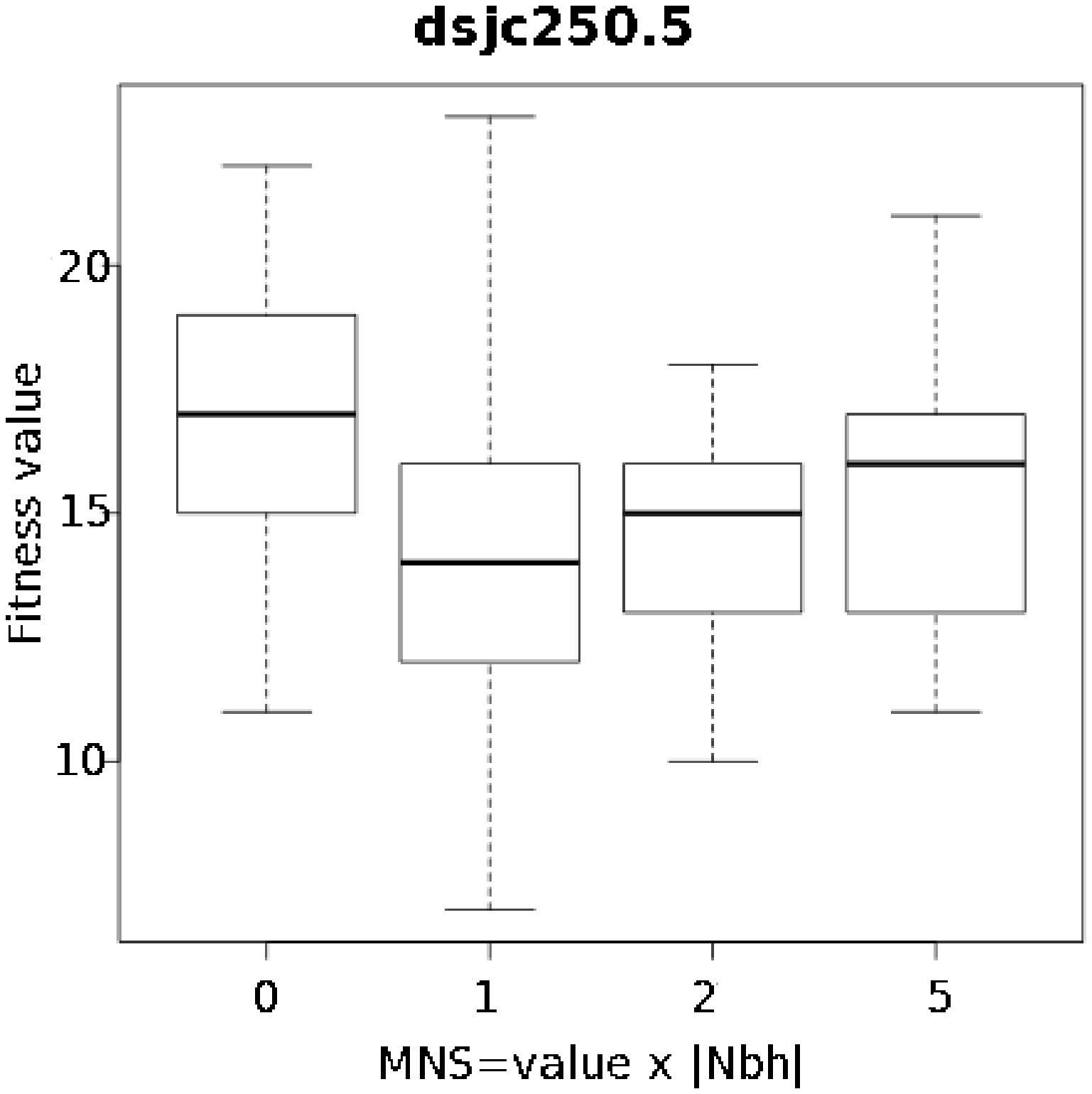} & \includegraphics[scale=0.35]{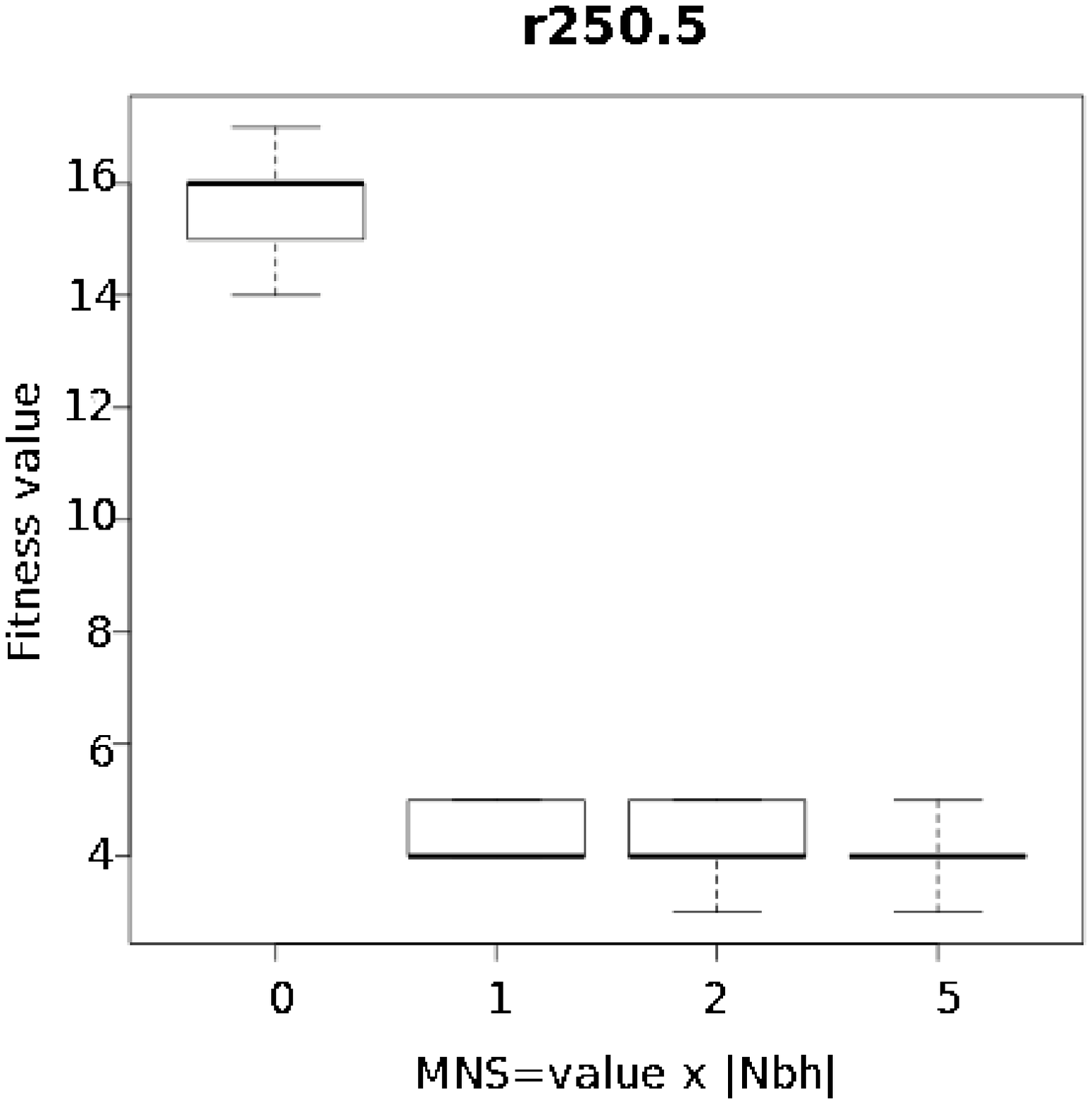} \\
  \includegraphics[scale=0.35]{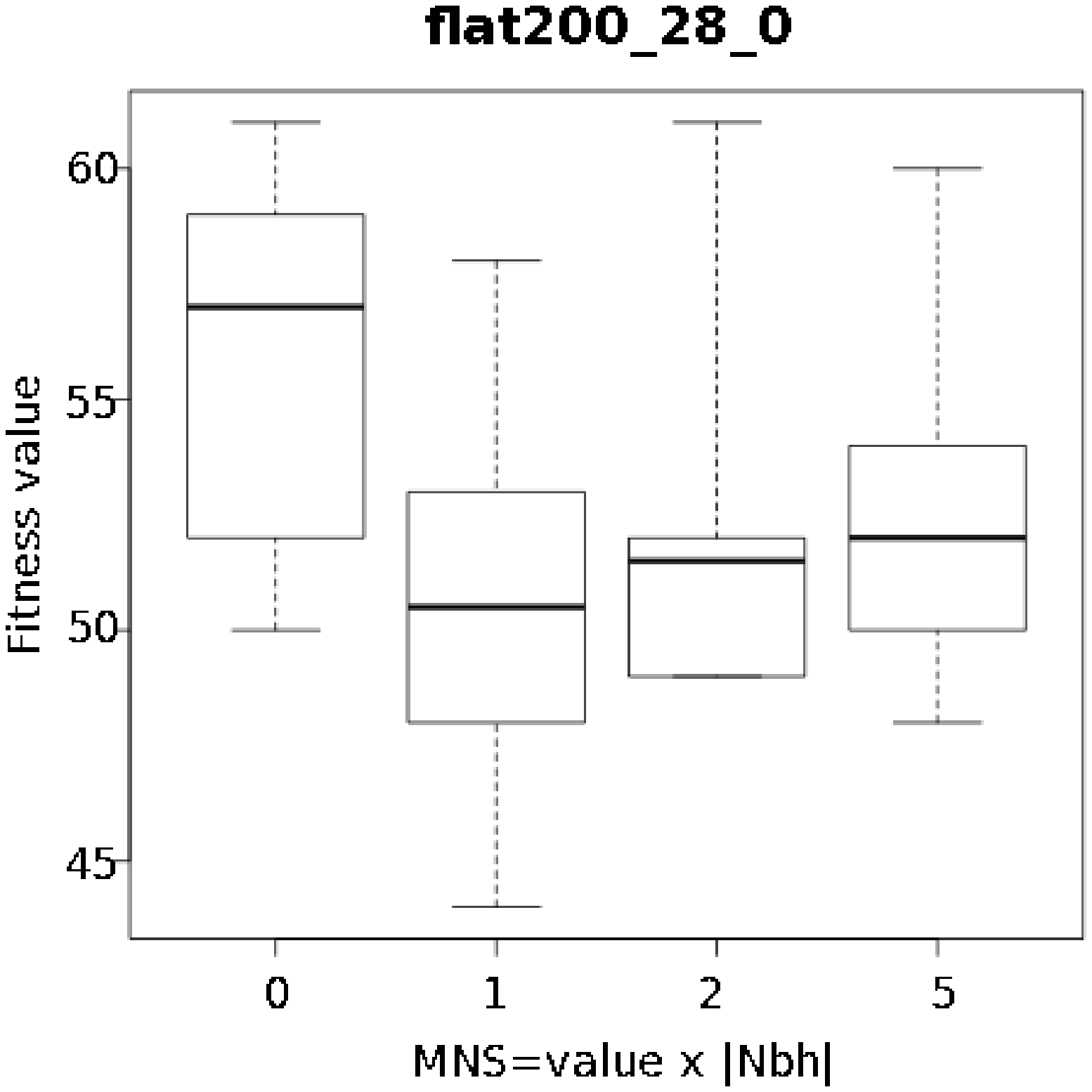} &  \includegraphics[scale=0.35]{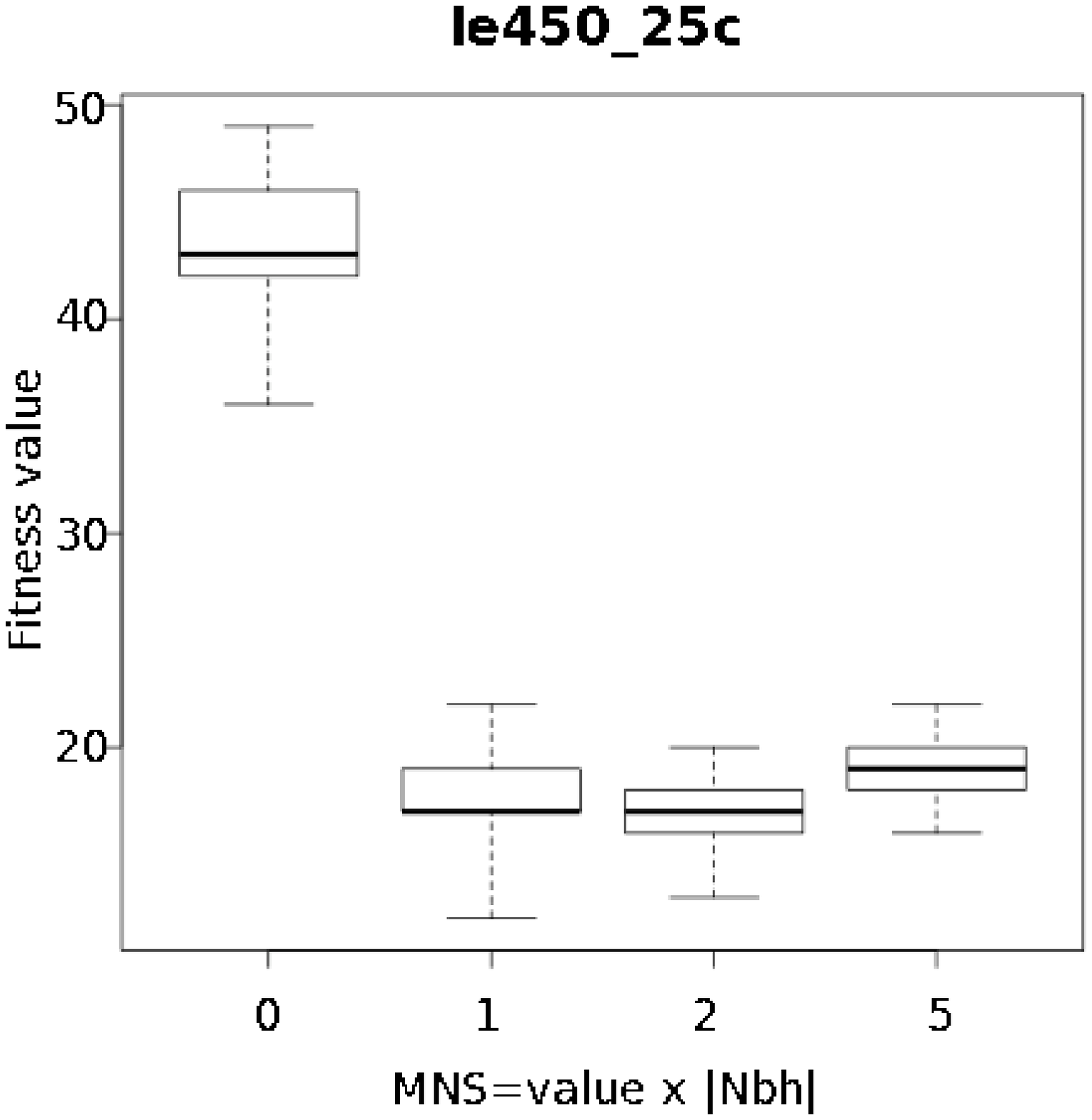} \\
  \end{tabular}
  \caption{ILS and NILS performance on 4 instances of the graph coloring problem.
  Next to the name of the instance is indicated the number of edges $m$ in the graph (as the fitness chosen
corresponds to the number of conflicting edges) and the chromatic value~$\chi$.}
  \label{GCP_boxplot_NILS_perf}
\end{figure}

Figure~\ref{GCP_boxplot_NILS_perf} presents the boxplot of the performance of the classical ILS ($MNS$ value equals to 0) and the three configurations of NILS ($MNS$ values equal to 1, 2 or 5 times the size of the neighborhood). Boxplots give the extent of the 30 fitness values found.
This figure shows first, that the performance of NILS are in average better than the ones of the classical ILS.
For the instances r250.5 and le\_450\_25c, the neutral degree ratios were high (respectively 9.7\% and 5.1\%), and the results are very promising 
as they show a clear improvement over the standard ILS. 
For the other instances, the neutral degree ratios were lower (1.2\% for dscj250.5, and 1.1\% for flat300\_28\_0), and the results obtained by NILS are only a little better than the classical ILS.
These results lead to the hypothesis that if the neutrality degree ratio of an instance is high, NILS will probably give good results. 
In other cases, it will not be attractive to use the neutrality, but however the results will not
be worse. 
Thus, exploiting the neutrality in the search process can lead to a better efficiency
and should not be discarded. \\
  
In these experiments, the performance of NILS is studied under different $MNS$ values for a same total number of evaluations. 
Results show that for $MNS$ values equal to 1 or 2 times the neighborhood size, performance is fairly similar. 
However, with a coefficient of 5, results are worse. 
That implies NILS can be stuck on plateaus on which searching portals is too expensive, 
and it may be preferable, in these cases, to escape the plateaus not to waste too much time. 

\subsection{Discussion}

The experimental results have shown that some instances of GCP present neutrality and that local search algorithms may be blocked on plateaus. 
Indeed, NILS with a $MNS$ equals to 0 is not able to find interesting solutions. 
However, when the neutrality is exploited in the local search, results are improved even if no configuration of NILS gives a legal solution.
This may be explained by the fact that these instances are the hardest instances of the literature, and for each, $k$ is set to the $\chi$-value, the best known chromatic number.
In 2010, Porumbel \textit{et al.}~\cite{porumbel:2010_2} made a comparison between their algorithm dedicated to GCP and the 10 best performing algorithms from the literature.
Except the Iterated Local Search, all the other algorithms are well-sophisticated and specific to GCP. Indeed, GCP-specific mechanisms are used to improve the search. These mechanisms require a huge knowledge on the GCP to be designed and tuned efficiently.
Despite this high level of sophistication, the comparison points out the difficulty for some algorithms to find the $\chi$-value.
For example, the results reported for the instances considered above indicate that:
The instance r250.5 is solved to the optimality ($k=\chi$) by only four algorithms out of six.
The instance le\_450\_25c is solved to the optimality only by six algorithms out of ten.
And, the instance flat300\_28\_0 is solved to the optimality only by  four algorithms out of eleven.
Moreover, the ILS~\cite{paquete:2002} never find the $\chi$-value for the two last instances.
Its performance illustrate the difficulty for a generic algorithm to be efficient. \\

This comparison highlights how these DIMACS instances are hard to solve when considering the $\chi$-value.
It may explain the performance of NILS described above.
The comparison between the ILS and the NILS performance on the DIMACS instances points out the interest of using neutrality.
However, one think that this neutrality should be exploited under the design of GCP-specific mechanisms.

\section{Conclusion}

This paper answers to the assumption that the GCP presents the neutrality property.
Indeed, for some hard instances of GCP, experiments show that a solution may have several neighbors with the same fitness value.
Then, the neutrality have to be taken into account when solving the GCP
as it may explain some good or bad performance of algorithms.
The neutrality characteristic is particularly interesting when it appears for local optima.
Local optima are integrated in a plateau of neighboring solutions with the same fitness value.
But, some of them are not local optima as they have at least one improving neighbor.
They are called portal.
Then, experiments points out the easiness to reach such portals from a local optima by moving randomly on the plateau.
From this observation, we perform NILS, a generic algorithm, exploiting the neutrality when is trapped on a local optimum.
It shows the benefit of taking the neutrality into account when solving these instances
as NILS gives better results than a classical ILS.\\

This paper should be considered as a preliminary work on the neutrality of the GCP.
Indeed, one points out the neutrality of some hard instances and gives the degree of this neutrality.
However, the performance of NILS are not as good as expected,
but, it shows the potential of exploiting neutrality to solve the GCP.
Since heuristic methods represent the state-of-the-art algorithms~\cite{caramia:2006,porumbel:2010},
one wants to investigate how to exploit neutrality in such heuristics.

\bibliographystyle{splncs03}
\bibliography{biblio,marmion}

\begin{thebibliography}{10}
\providecommand{\url}[1]{\texttt{#1}}
\providecommand{\urlprefix}{URL }

\bibitem{avanthay:2003}
Avanthay, C., Hertz, A., Zufferey, N.: A variable neighborhood search for graph
  coloring. European Journal of Operational Research  151(2),  379--388 (2003)

\bibitem{bastolla:2003}
Bastolla, U., Porto, M., Roman, H.E., Vendruscolo, M.: Statistical properties
  of neutral evolution. Journal Molecular Evolution  57(S),  103--119 (2003)

\bibitem{bouziri:2011}
Bouziri, H., Mellouli, K., Talbi, E.G.: The k-coloring fitness landscape.
  Journal of Combinatorial Optimization  21(3),  306--329 (2011)

\bibitem{brelaz:1979}
Br{\'{e}}laz, D.: {New methods to color the vertices of a graph}.
  Communications of the ACM  22(4),  251--256 (1979)

\bibitem{briggs:1994}
Briggs, P., Cooper, K.D., Torczon, L.: Improvements to graph coloring register
  allocation. ACM Transactions on Programming Languages and Systems  16,
  428--455 (1994)

\bibitem{caramia:2006}
Caramia, M., Dell'Olmo, P., Italiano, G.F.: Checkcol: Improved local search for
  graph coloring. Journal of Discrete Algorithms pp. 277--298 (2006)

\bibitem{chams:1987}
Chams, M., Hertz, A., de~Werra, D.: Some experiments with simulated annealing
  for coloring graphs. European Journal of Operational Research  32(2),  260 --
  266 (1987)

\bibitem{fleurent:1996}
Fleurent, C., Ferland, J.: Genetic and hybrid algorithms for graph coloring.
  Annals of Operations Research  63,  437--461 (1996)

\bibitem{galinier:2006}
Galinier, P., Hertz, A.: A survey of local search methods for graph coloring.
  Computers and Operations Research  33(9),  2547--2562 (2006)

\bibitem{garey:1990}
Garey, M.R., Johnson, D.S.: Computers and Intractability; A Guide to the Theory
  of NP-Completeness. W. H. Freeman \& Co. (1990)

\bibitem{hertz:1987}
Hertz, A., de~Werra, D.: Using tabu search techniques for graph coloring.
  Computing  39(4),  345--351 (1987)

\bibitem{maniezzo:2000}
Maniezzo, V., Carbonaro, A.: An ants heuristic for the frequency assignment
  problem. Future Generation Computer Systems  16(8),  927 -- 935 (2000)

\bibitem{marmion@evocop:2011}
Marmion, M.É., Dhaenens, C., Jourdan, L., Liefooghe, A., Verel, S.: {NILS}: a
  neutrality-based iterated local search and its application to flowshop
  scheduling. In: Proceedings of the 12th European Conference of Evolutionary
  Computation in Combinatorial Optimization. pp. 191-- 202. EvoCOP 2011, LNCS,
  Springer (2011)

\bibitem{marmion@lion:2011}
Marmion, M.É., Dhaenens, C., Jourdan, L., Liefooghe, A., Verel, S.: On the
  neutrality of flowshop scheduling fitness landscapes. In: Proceedings of the
  5th Learning and Intelligent OptimizatioN Conference. pp. 238--252. LION
  2011, LNCS, Springer (2011)

\bibitem{marmion@jmma:2011}
Marmion, M.É., Jourdan, L., Dhaenens, C.: Fitness landscape analysis and
  metaheuristics efficiency. Journal of Mathematical Modelling and Algorithms
  (2013, to appear)

\bibitem{merz:1997}
Merz, P., Freisleben, B.: Memetic algorithms for the traveling salesman
  problem. Complex Systems  13,  297--345 (1997)

\bibitem{merz:2000_2}
Merz, P., Freisleben, B.: Fitness landscape analysis and memetic algorithms for
  the quadratic assignment problem. IEEE Transactions on Evolutionary
  Computation  4(4),  337--352 (2000)

\bibitem{paquete:2002}
Paquete, L., Stutzle, T.: An experimental investigation of iterated local
  search for coloring graphs. In: Applications of Evolutionary Computing,
  volume 2270 of LNCS. pp. 122--131. Springer-Verlag (2002)

\bibitem{porumbel:2010_2}
Porumbel, D.C., Hao, J.K., Kuntz, P.: An evolutionary approach with diversity
  guarantee and well-informed grouping recombination for graph coloring.
  Computers and Operations Research  37(10),  1822--1832 (2010)

\bibitem{porumbel:2010}
Porumbel, D.C., Hao, J.K., Kuntz, P.: A search space "cartography" for guiding
  graph coloring heuristics. Computers and Operations Research  37(4),
  769--778 (2010)

\bibitem{reidys:2001}
Reidys, C.M., Stadler, P.F.: Neutrality in fitness landscapes. Applied
  Mathematics and Computation  117,  321--350 (2001)

\bibitem{stadler:1996}
Stadler, P.F.: Landscapes and their correlation functions. Journal of
  Mathematical Chemistry  20,  1--45 (1996)

\bibitem{tricjk:2007}
Trick, M., Yildiz, H.: A large neighborhood search heuristic for graph
  coloring. In: Integration of AI and OR Techniques in Constraint Programming
  for Combinatorial Optimization Problems, LNCS, vol. 4510, pp. 346--360.
  Springer Berlin / Heidelberg (2007)

\bibitem{verel:2004}
{V}erel, S., {C}ollard, P., {C}lergue, M.: {S}cuba {S}earch : when selection
  meets innovation. In: Proceedings of the 2004 Congress on Evolutionary
  Computation. pp. 924 -- 931. CEC'04, {IEEE} {P}ress (2004)

\bibitem{verel:2007}
{V}erel, S., Collard, P., Tomassini, M., Vanneschi, L.: Fitness landscape of
  the cellular automata majority problem: view from the ``{O}lympus''.
  Theoretical Computer Science  378,  54--77 (2007)

\bibitem{weinberger:1990}
Weinberger, E.: Correlated and uncorrelated fitness landscapes and how to tell
  the difference. Biological cybernetics  63,  325--336 (1990)

\bibitem{wilke:2001}
Wilke, C.O.: Adaptative evolution on neutral networks. Bulletin of Mathematical
  Biology  63,  715--730 (2001)

\bibitem{wright:1932}
Wright, S.: The roles of mutation, inbreeding, crossbreeding and selection in
  evolution. In: Proceedings of the Sixth International Congress on Genetics.
  vol.~1 (1932)

\bibitem{zufferey:2008}
Zufferey, N., Amstutz, P., Giaccari, P.: Graph colouring approaches for a
  satellite range scheduling problem. J. Scheduling  11(4),  263--277 (2008)

\end{thebibliography}

\end{document}